\documentclass[11pt]{article}

\usepackage{amssymb,amsmath,amsthm}
\usepackage{graphicx, tikz}

\usepackage{color}
\definecolor{darkgreen}{rgb}{0,.4,0.2}
\definecolor{darkagenta}{rgb}{.5,0,.5}
\definecolor{darkred}{rgb}{0.85,0,0}
\definecolor{darkblue}{rgb}{0,0,.6}
\definecolor{lightgray}{gray}{.95}
\definecolor{rgrey}{rgb}{.8,0.4,.4}  
\definecolor{grey}{rgb}{.13,.13,.13}  

\newcommand{\Id}{ {\rm Id}}

\newcommand{\xF}{x_{\rm \scriptstyle F}}

\newcommand{\triple}[3]{\fbox{$\begin{aligned} w &  #1\\[-.1cm] C_w &  #2\\[-.1cm] R_w &  #3  \end{aligned}$}}

\newtheorem*{theorem*}{Theorem}
\newtheorem{theorem}{Theorem}[section]

\newtheorem{lemma}[theorem]{Lemma}

\newtheorem{proposition}[theorem]{Proposition}

\advance \topmargin by -\headheight
\advance \topmargin by -\headsep
\evensidemargin \oddsidemargin
\begin{document}

\begin{centering}

{\huge \textbf{The structure of Zeckendorf expansions}}

\bigskip

{\bf \large F.~Michel Dekking}

\bigskip

{ DIAM,  Delft University of Technology, Faculty EEMCS,\\ P.O.~Box 5031, 2600 GA Delft, The Netherlands.}

\medskip

{\footnotesize \it Email:  F.M.Dekking@TUDelft.nl}

\end{centering}

\medskip

\center{\small June 12, 2020}

\begin{abstract}
 \noindent  In  this paper we classify the Zeckendorf expansions according to their digit blocks. It turns out that if we consider these digit blocks as labels on the Fibonacci tree, then the  numbers ending with a given digit block in their Zeckendorf expansion appear as compound Wythoff sequences in a natural way on this tree. Here the digit blocks consisting of only $0$'s are an exception. We also give a second description of these occurrence sequences as  generalized Beatty sequences. Finally, we characterize the numbers with a fixed digit block occurring at an arbitrary fixed position in their Zeckendorf expansions, and determine their densities.
\end{abstract}

\medskip


\bigskip

\section{Introduction}

We define the Zeckendorf expansion, as introduced in  \cite{Lekker} and \cite{Zeck}. Let $F_0=0,\, F_1=1,\, F_2=1,\dots$ be the Fibonacci numbers.
Let $\ddot{F}_0=1,\, \ddot{F}_1=2, \dots$ be the twice shifted Fibonacci numbers, defined by  $\ddot{F}_i=F_{i+2}$.
  Ignoring leading and trailing zeros, any  natural number $N$ can be written uniquely as
  $$N= \sum_{i=0}^{\infty} d_i \ddot{F}_i,\vspace*{-.0cm}$$
  with digits $d_i=0$ or 1, and where $d_id_{i+1} = 11$ is not allowed.
We denote the Zeckendorf expansion of $N$ as $Z(N)$, writing $Z(N)=\dots d_2d_1d_0$.

  \medskip

%
%
%

In Theorem \ref{th:Zeck-w} in Section \ref{sec:main} we characterize those numbers $N$ which have a Zeckendorf expansion ending with the digit block $w=d_{m-1}\dots d_1d_0$ for any $m$, and any choice of digits.

Several authors have obtained special cases of our results. In, e.g., the paper \cite{Rytter}, the digit blocks $w=0^m$, for $m\ge 1$ are treated.
It is interesting (in the light of our Proposition \ref{prop:A-B}), that Rytter (\cite{Rytter}, page 219) writes ``The remarkable property of the sequence of occurrences of a word $w=0^m$ is that its difference sequence is structurally isomorphic to the infinite Fibonacci word\dots".

\smallskip

In Section \ref{sec:gen-pos} we solve the general problem of characterizing those numbers $N$ which have a Zeckendorf expansion $Z(N)=\dots d_k\dots d_2d_1d_0$ such that the digit block $ww_{m-1}\dots w_0$ occurs at position $k$ in $Z(N)$, i.e., $d_{k+m-1}\dots d_k=w_{m-1}\dots w_0$.

\smallskip

In Section \ref{sec:dens} we given the densities $f^{(k)}_w$ of the numbers $N$ which have a Zeckendorf expansion $Z(N)=\dots d_k\dots d_2d_1d_0$ such that the digit block $w$ occurs at position $k$.

Here too,  a special case has been considered in the literature: in the paper \cite{Griffith-2010} the densities $f^{(k)}_1$ of the word $w=1$ in position $k$ are computed (in a different way).

\bigskip

The question arises what one can say about more general numeration systems, for example, about Ostrowski representations of the natural numbers.
It appears that even for quadratic irrationals, the situation is much more complicated than for the golden mean.

One reason is that compound Beatty sequences, in particular iterated Beatty sequences, will no longer be generalized Beatty sequences in general. The well known example is $AA(n)=\lfloor\lfloor n\alpha \rfloor\alpha \rfloor$, with $\alpha=\sqrt{2}$, see, e.g., Theorem 1 in  \cite{Dekking-TCS}. So a generalization of Theorem \ref{thm:Car} does not exist, nor a `compound Beatty part' of our main Theorem \ref{th:Zeck-w}.

Nevertheless, expressions involving generalized Beatty sequences are possible (see, e.g., Theorem 5.3  in \cite{car-Pell}, which, as Theorem 7 and 8 in \cite{car-sco-hog}, is a special case of the problem for arbitrary digit blocks.).  These  are the subject of future work.

\section{Zeckendorf expansions}\label{sec:Zeck}

 Let the golden mean be given by  $\varphi:=(1+\sqrt{5})/2$. It is well-known (see, e.g., \cite{car-sco-hog}, \cite{Kimberling-1983}) that the numbers $N$ whose Zeckendorf expansion $Z(N)=\dots d_2d_1d_0$ has digit $d_0=0$ are exactly the elements $N=0,2,3,5,7,8,10,\dots$  from the lower Wythoff sequence
 \begin{equation}\label{eq:A}
 (A(N))=(\lfloor N\varphi \rfloor)=(1,3,4,6,8,9,11,\dots),
 \end{equation}
 with  1 subtracted.  Those with digit $d_0=1$ are exactly the elements $N=1,4,6,9,12,14,\dots$ from the upper Wythoff sequence
 \begin{equation}\label{eq:B}
 (B(N))=(\lfloor N\varphi^2 \rfloor)=(2,5,7,10,13,15,\dots),
 \end{equation}
 with  1 subtracted.

 We consider the general question: given any word $w$ of length $m$, what are the numbers $N$ whose Zeckendorf expansion $Z(N)=\dots d_2d_1d_0$ ends with the digits $d_{m-1}\dots d_0=w$?

\subsection{Compound Wythoff sequences}

An important role is played by  compositions of the two sequences $A$ and $B$ in Equation (\ref{eq:A}) and (\ref{eq:B}), also known as compound Wythoff sequences.
As usual, we write these compositions as words over the monoid generated by $A,B$. For example, the compound sequence  $AB$ is given by
$AB(N)=A(B(N))$ for $N=1,2\dots.$ These compound Wythoff sequences have been extensively studied, as, e.g., in \cite{car-sco-hog} and \cite{fraenkel1994}.

It turns out that, with  exception of the words $w=0^m$, the numbers $N$ whose Zeckendorf expansion $Z(N)$ ends with the digits $d_{m-1}\dots d_0=w$,  are given by a compound Wythoff sequence, which we denote by $C_w$ or $C(w)$.

\smallskip

In the pioneering paper \cite{car-sco-hog} by Carlitz, Scoville  and Hoggatt, we find that for $m\ge 0$
\begin{align}
C(10^{2m+1})  & =B^{m+1}A,\;\;\;\quad C(10^{2m})=AB^mA,\label{eq:WA}\\
C(0010^{2m+1})& =B^{m+1}AA, \;C(010^{2m})=AB^mAA,\label{eq:w0010}\\
C(1010^{2m+1})& =B^{m+1}AB,\;C(1010^{2m})=AB^mAB.\label{eq:w1010}
\end{align}
These are given in their Theorems 7 and 8.
It is remarkable that these results are presented as their main results in their introduction, but that their Theorem 13 (see Theorem \ref{thm:Car} below), which we consider the most important result in \cite{car-sco-hog},  is not mentioned.

The successive  compound Wythoff sequences  $AA, BA, ABA, BBA,\dots$ in Equation (\ref{eq:WA}) are the successive columns of the so-called Wythoff array $W$, defined by
$$W(n,m)=F_{m+1}\lfloor n\varphi \rfloor + (n -1)F_{m} \qquad n\ge 1, m \ge 0.$$
This is stated  in Theorem 10 in the paper \cite{Kimberling-JIS}, but  already implicitly proved in \cite{car-sco-hog}, see Equation (3.8), (6.10) and (6.16) in that paper. We remark further that Proposition \ref{prop:A-B} gives a very simple way to obtain  (\ref{eq:w0010}) and (\ref{eq:w1010}) from (\ref{eq:WA}).

\smallskip

\begin{theorem}{\bf (Carlitz-Scoville-Hoggatt)}\label{thm:Car}
Let $U = (U(n))_{n \geq 1}$ be a composition of the Wythoff sequences
$A$ and $B$, containing $i$ occurrences of $A$ and $j$ occurrences of $B$, then
$$
U(n) = F_{i + 2j}\,A(n) + F_{i+2j-1}\,n - \lambda_U,\quad \text{for all }  n \geq 1,
$$
where $F_k$ are the Fibonacci numbers and $\lambda_{U}$ is a constant.
\end{theorem}

\subsection{Generalized Beatty sequences}

Let $\alpha$ be an irrational number larger than 1. We call any sequence $V$ with terms of the form $V(n) = p\lfloor n \alpha \rfloor + q n +r $, $n\ge 1$ a \emph{generalized Beatty sequence}. Here $p,q$ and $r$ are integers, called the \emph{parameters} of $V$.

Theorem  \ref{thm:Car} has a useful extension, given in the paper \cite{GBS}. In its statement below, as Lemma \ref{corr:VA}, a typo in its source is corrected.

\begin{lemma}{\bf (\cite{GBS}, Corollary 2) }\label{corr:VA} Let $V$ be a generalized Beatty sequence with parameters $(p,q,r)$, and $\alpha=\varphi$. Then $VA$ and $VB$ are generalized Beatty sequences with parameters $(p_{V\!A},q_{V\!A},r_{V\!A})=(p+q,p,r-p)$ and $(p_{V\!B},q_{V\!B},r_{V\!B})=(2p+q,p+q,r)$.
\end{lemma}

\subsection{Main theorem}\label{sec:main}

To formulate our result for general digit blocks $w$, it is convenient to add $0$'s to the expansion of a number $N$ in $\{0,\dots,F_{n}-1\}$ such that the total length of the word $Z(N)$ becomes $n-2$. We denote this word as $Z^*(N)$. For example, for $n=5$, we have
\begin{align*} Z(0)&=0, Z(1)=1, Z(2)=10, Z(3)=100, Z(4)=101,\\
Z^*(0)&=000, Z^*(1)=001, Z^*(2)=010, Z^*(3)=100, Z^*(4)=101.
\end{align*}
Note that $Z^*(\cdot)$ depends on the particular $n$ that one considers, but we do not add this to avoid burdening the notation. The appropriate $n$ will always be clear from the context.

In the following theorem, occurrences of a word $w$ have to be interpreted in the $Z^*$-sense.

\medskip

\begin{theorem}\label{th:Zeck-w}
\noindent For any natural number $m$ fix a word $w$ of $0$'s and $1$'s, containing no $11$. Then---except if $w=1$, or $w=0^m$---the sequence $R_w$ of occurrences of numbers $N$ such that  the $m$ lowest digits of the Zeckendorf expansion of $N$ are equal to $w$, i.e., $d_{m-1}\dots d_0=w$, is a compound Wythoff sequence $C_w$.\\
In the exceptional cases  $w=1$, we have  $R_w=B-1$;  when $w=0^m$, we have  $R_w=A^m-1$.\\
The representation of $R_w$ as a generalized Beatty sequence is given, without exception, by\\[-.4cm]
$$R_w=F_{m}A+F_{m-1}\Id+\gamma_w  \quad \text{or by} \quad R_w=F_{m+1}A+F_m\Id+\gamma_w,$$\\[-.6cm]
for some negative  integer $\gamma_w$.
The first representation holds for all $w$ starting with $w_{m-1}=0$, the second for all $w$ with $w_{m-1}=1$.
\end{theorem}

\noindent Theorem \ref{th:Zeck-w} will be proved in Section \ref{sec:proof}.

\smallskip

 The words $w$ without $11$ are naturally ordered in a tree, the Fibonacci tree\footnote{These are different from the `Fibonacci trees' considered in \cite{Capocelli}}.
The first four levels of this tree are depicted below. The nodes are labeled with the $w$'s, the corresponding compound Wythoff sequences $C_w$, and the $R_w$'s, expressed as generalized Beatty sequences.

\bigskip

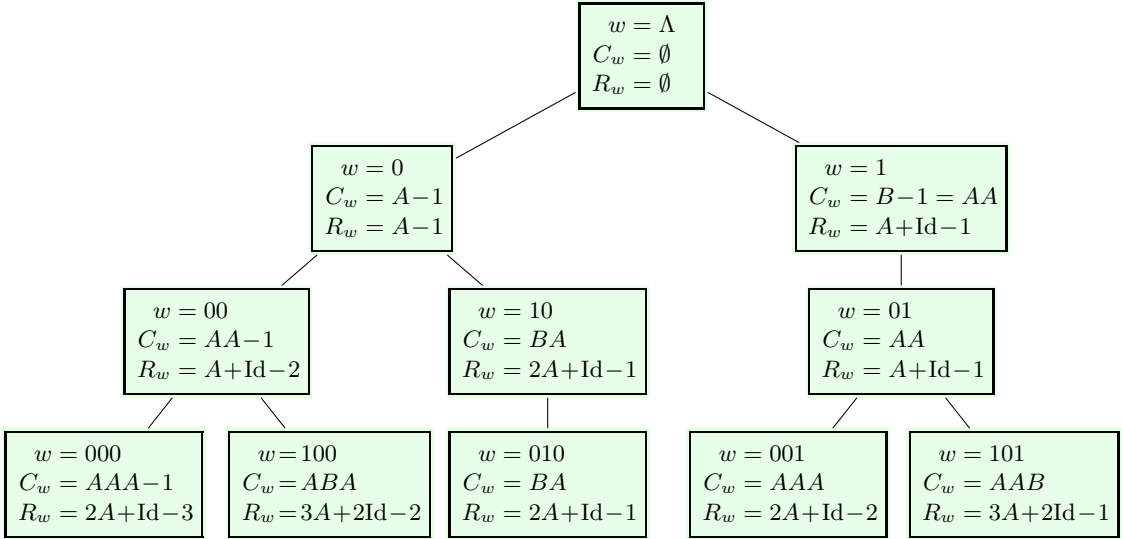
\begin{tikzpicture}
[level distance=19mm,
every node/.style={fill=green!10,rectangle,inner sep=1pt},
level 1/.style={sibling distance=69mm,nodes={fill=green!10}},
level 2/.style={sibling distance=44mm,nodes={fill=green!10}},
level 3/.style={sibling distance=30mm,nodes={fill=green!10}}]
\node {\footnotesize  $\triple{=\Lambda}{=\emptyset\quad}{=\emptyset\quad}$}            
child {node {\footnotesize $\triple{=0}{=A\!-\!1}{=A\!-\!1}$}            
 child {node {\footnotesize $\triple{=00}{=AA\!-\!1}{=A\!+\!\Id\!-\!2}$}          
  child {node {\footnotesize $\triple{=000}{=AAA\!-\!1}{=2A\!+\!\Id\!-\!3}$}}       
  child {node {\footnotesize $\triple{\!=\!100}{\!=\!ABA}{\!=\!3A\!+\!2\Id\!-\!2}$}}    
}
 child {node {\footnotesize $\triple{=10}{=BA}{=2A\!+\!\Id\!-\!1}$}         
  child {node {\footnotesize $\triple{=010}{=BA}{=2A\!+\!\Id\!-\!1}$}}     
}
}
child {node {\footnotesize $\triple{=1}{=B\!-\!1=AA}{=A\!+\!\Id\!-\!1}$}            
 child {node {\footnotesize $\triple{=01}{=AA}{=A\!+\!\Id\!-\!1}$}          
  child {node {\footnotesize $\triple{=001}{=AAA}{=2A\!+\!\Id\!-\!2}$}}     
   child {node {\footnotesize $\triple{=101}{=AAB}{=3A\!+\!2\Id\!-\!1}$}}  
}
};
\end{tikzpicture}

\bigskip

\subsection{The basic recursion}

We partition the natural numbers in sets $\Lambda_n$, and consider sets $\Psi_n$ given by
$$\Lambda_n :=\{F_n,\dots,F_{n+1}\!-\!1\},\quad \Psi_n :=\{0,\dots,F_{n}\!-\!1\},\quad n\ge 2.$$
Note that the elements from $\Lambda_n$ are exactly the numbers with $n-1$ digits in their Zeckendorf expansions, and that the elements of $\Psi_n$ are the numbers with $n-2$ digits, or less, in their Zeckendorf expansions.

 \smallskip

Immediately from the definition of $Z$-expansions one obtains that the following basic recursion holds.

\begin{lemma}\label{lem:recurs}
If $N\in \Lambda_n$, then $Z(N)=1\,Z^*(N\!-F_n)$.
\end{lemma}

\noindent For example, $F_6=8$, $\Lambda_6=\{8,9,10,11,12\}$, and $Z(11)=10100=1Z^*(3)$.

\subsection{The role of the infinite Fibonacci word}

Note that we index the elements of the digit blocks $w$ in reverse order, to comply with the order of the digits in  $Z(N)$, and with the order of the levels of the Fibonacci tree.

Let $w=w_{m-1}\dots w_0$  be a word with  $w_{m-1}=0$. The idea is to determine how the  occurrences of the numbers $N$ with $Z(N)=\dots 0w$, and those with  $Z(N)=\dots 1w$,  are intertwined.

\medskip

Let $f$ be the Fibonacci morphism on the alphabet $\{a,b\}$ given by $f(a) = ab, \; f(b) =a.$ This morphism generates an infinite word $\xF=abaababaaba\dots$ by iteration, see e.g., the monograph \cite{Lothaire}, Proposition 1.2.8 and Example 1.2.10.

\smallskip

\begin{proposition}\label{prop:fib} For any natural number $m>1$ fix a word $w=w_{m-1}\dots w_0$ of $0$'s and $1$'s, containing no $11$, with $w_{m-1}=0$. Code any occurrence of $0w$ at the end of a $Z(N)$ by $a$, and any occurrence of $1w$ at the end of a $Z(N)$ by $b$, in the order of these occurrences in $\Psi_{m+n}$, then the resulting word is equal to $f^{n-2}(a)$, for $n=3,4,\dots$.
\end{proposition}

\noindent {\it Proof:} We prove this by induction on $n$.

\smallskip

Because the length of $0w$ and $1w$ is  $m+1$, there is a  single occurrence of $0w$ in  $\Psi_{m+3}$,  and a single occurrence of $1w$ in $\Psi_{m+3}$, in that order. The corresponding coding is $ab=f(a)$. This deals with the case $n=3$.

\smallskip

All the numbers $N$ in $\Psi_{m+4}$, but not in $\Psi_{m+3}$, have an expansion $Z(N)$ starting with 1. This means that from the three possible left extensions of $0w$ and $1w$, only $10w$ will occur as a $Z(N)$, with $N$ in $\Psi_{m+4}$, obviously after the occurrences of $0w$ and $1w$. So the coding of the occurrences in $\Psi_{m+4}$ equals $aba=f^2(a)$. This was the case $n=4$.

\smallskip

Now consider $\Psi_{m+n}$ for $n\ge5$. By Lemma \ref{lem:recurs} this set can be written as
$$\Psi_{m+n}= \Psi_{m+n-1} \cup \big(\Psi_{m+n-2} + F_{m+n-1}\big).$$
Here we write for a set $E$ and a number $x$, $E+x:=\{e+x:e\in E\}$.

Note that the occurrences of $0w$ and $1w$ as end blocks of Zeckendorf expansions of numbers in $\Psi_{m+n-2} + F_{m+n-1}$ are the same as for those in $\Psi_{m+n-2}$. From the induction hypothesis it then follows that the sequence of occurrences of $0w$ and $1w$ as end blocks is coded by the word
$$f^{n-3}(a)f^{n-4}(a)=f^{n-3}(a)f^{n-3}(b)=f^{n-3}(ab)=f^{n-2}(a).$$
This finishes the induction proof. \hfill $\Box$

\bigskip

The main part of Theorem \ref{th:Zeck-w} is a consequence of the following result.

\begin{proposition}\label{prop:A-B} For any natural number $m>1$ fix a word $w=w_{m-1}\dots w_0$ of $0$'s and $1$'s, containing no $11$, with $w_{m-1}=0$. Let $C_w$ be the Wythoff-coding of the sequence of occurrences of the numbers $N$ whose $Z^*$-expansion ends with  $w$. Then $C_{0w}=C_w A$, and $C_{1w}=C_w B$.
\end{proposition}

\noindent  {\it Proof:} One recalls (see, e.g., \cite{Lothaire}) that the letters $a$ in the infinite Fibonacci word $\xF$, which has the $f^n(a)$ as prefix, occur at positions given by the lower Wythoff sequence $A$, and the letters $b$ occur at positions given by $B$. Now the proposition follows directly from Proposition \ref{prop:fib}. \hfill $\Box$

\subsection{Two particular cases of digit blocks}

The digit block $w=0^m$ behaves exceptionally in Theorem \ref{th:Zeck-w}, and also the digit block $w=10^m$ needs special care.

\begin{lemma}\label{lem:0m}
Let $A$ be the lower Wythoff sequence, and $B$ the upper Wythoff sequence. Then
$$(A^m\!-\!1) A = A^{m+1}\!-\!1,\;\; (A^{2m\!-\!1}\!-\!1) B = B^m A,\;\; (A^{2m}\!-\!1) B = AB^m A, \quad \text{for } m\ge 1.$$
\end{lemma}

\noindent {\it Proof:} From Theorem \ref{thm:Car}, filling in $A^m(1)=1$, we obtain for all $m$
$$A^m(n)\!-\!1=F_m\,A(n)+F_{m-1}\,n-F_{m+1}.$$
Applying Lemma \ref{corr:VA}, this yields
$$A^m(A(n)\!-\!1)=(F_m+F_{m-1})\,A(n)+F_m\,n-F_m-F_{m+1}=F_{m+1}\,A(n)+F_m\,n-F_{m+2}.$$
This is indeed equal to $A^{m+1}\!-\!1$.

Applying Lemma \ref{corr:VA}, now to $(A^{2m\!-\!1}\!-\!1) B$, we obtain
$$A^{2m\!-\!1}(B(n))-1 = (2F_{2m-1}+F_{2m-2})\,A(n)+(F_{2m-1}+F_{2m-2})\,n-F_{2m}=F_{2m+1}\,A(n)+F_{2m}n-F_{2m}.$$
On the other hand, we find with Theorem \ref{thm:Car}, and by using Lemma \ref{corr:VA} appropriately, that
$$B^m(A(n)) = F_{2m+1}\,A(n)+F_{2m}n-\lambda_{B^mA}= F_{2m+1}\,A(n)+F_{2m}\,n-F_{2m}.$$
This establishes the second equation.

For the third we compute
$$A(B^m(A(n))=F_{2m+2}A(n)+F_{2m+1}n-\lambda_{AB^mA}= F_{2m+2}A(n)+F_{2m+1}n-F_{2m+1},$$
which establishes the third equation.   \hfill $\Box$

\subsection{Proof of Theorem \ref{th:Zeck-w}}\label{sec:proof}

The proof is by induction on the length $m$ of the digit block $w$. For $m=1$, we have $R_0=A-1$ and $R_1=A+\Id-1$, since $B=A+\Id$. Cf.~Equations (\ref{eq:A}) and (\ref{eq:B}).

Suppose we know the result for all words $w$ of length $m$. We prove that it holds for all words $w$ of length $m+1$. It holds for $w=0^{m+1}$, by Proposition \ref{prop:A-B}, Lemma \ref{lem:0m}, and Lemma \ref{corr:VA}. For $w=10^{m}$ we also apply these three results, distinguishing between odd and even $m$.
For the other words $w$ of length $m$, we distinguish between  $w=0u$ or $w=1u$, where $u$ is a word of length $m-1$.

The word $1u$ just generates the single word $w=01u$, which has $C_w=C_{1u}$,
and $R_w=R_{1u}=F_{m+1}A+F_m\Id+\gamma_{1u}$, which is correct as announced, since $w$ starts with 0.

The word $0u$  generates the two words $w=00u$ and $10u$, which by Proposition \ref{prop:A-B} are compound Wythoff. Also,  the induction hypothesis is that $R_{0u}=F_{m}A+F_{m-1}\Id+\gamma_{0u}$. By  Lemma \ref{lem:0m} this implies that $R_{00u}= F_{m+1}A+F_{m}\Id+\gamma_{00u}  $ and
$R_{10u}= F_{m+2}A+F_{m+1}\Id+\gamma_{10u}$. This is exactly what had to be proved.  \hfill $\Box$

\medskip

 The value of $\gamma_w$ in the representation $R_w=F_{m}A+F_{m-1}\Id+\gamma_w$  or $R_w=F_{m+1}A+F_m\Id+\gamma_w$ in Theorem \ref{th:Zeck-w} can be easily computed.

 \begin{proposition}\label{prop:gamma} For any natural number $m>1$ fix a word $w=w_{m-1}\dots w_0$ of $0$'s and $1$'s, containing no $11$.
 Let $T_{00}:=\{0< k<m: w_kw_{k-1}=00\}.$ Then
$$-\gamma_w=1+\sum_{k\in T_{00}}F_k.$$
\end{proposition}

\noindent {\it Proof:} This is implied directly by $\gamma_0=\gamma_1=-1$, and Lemma \ref{corr:VA}, following the steps in the proof of Theorem \ref{th:Zeck-w}.  \hfill $\Box$

\subsection{Digit blocks at arbitrary positions}\label{sec:gen-pos}

 The general question is: what is the sequence $R^{(k)}_w$ of occurrences of numbers $N$ such that  the length $m$ digit block $w$ ends at position $k$ in  the Zeckendorf expansion of $N$? For $k=0$ the answer is given by Theorem \ref{th:Zeck-w}.

When $(a_n)$ and $(b_n)$ are two increasing sequences, indexed by $\mathbb{N}$, then we mean by the union of $(a_n)$ and $(b_n)$ the increasing sequence whose terms go through the set $\{a_n, b_n: n\in \mathbb{N}\}$. By iteration, we also consider arbitrary finite unions of increasing sequences.

\medskip

\begin{theorem}\label{th:gen-pos}
\noindent For any natural number $m$ fix a word $w$ of $0$'s and $1$'s, containing no $11$. Let k be a positive  integer.  Let $v$ be the word $v=w\,0^k$. Then the sequence $R^{(k)}_w$ of occurrences of numbers $N$ with expansion $Z(N)=\dots d_k\dots d_1d_0$ such that   $d_{k+m-1}\dots d_k=w=w_{m-1}\dots w_0$, is a union  of $F_{k+2-w_0}$ generalized Beatty sequences, given by
$$ F_{k+m+w_{m-1}}A+F_{k+m-1+w_{m-1}}\Id+\gamma_v,\,\dots,\, F_{k+m+1}A+F_{k+m}\Id+\gamma_v+F_{k+2-w_0}-1.$$
\end{theorem}

\noindent {\it Proof:} There are $F_{k+2}$ words of length $k$, containing no $11$.  So this theorem is implied directly by  Theorem \ref{th:Zeck-w}: a digit block  $w=d_{k+m-1}\dots d_{k+1}\,0$ extends to $F_{k+2}$ digit blocks $d_{k+m-1}\dots d_{k+1}\,0\,d_{k-1}\dots d_0$, whereas a digit block  $w=d_{k+m-1}\dots d_{k+1}\,1$ extends to $F_{k+1}$ digit blocks $d_{k+m-1}\dots d_{k+1}\,10\,d_{k-2}\dots d_0$. The corresponding $N$'s are consecutive, with the smallest such $N$ equal to $v=w\,0^k$.
\hfill $\Box$

\medskip

\noindent  {\bf Example} Let $w=00$, and $k=2$. Then $d_3d_2=00$ if and only if $d_3d_2d_1d_0=0000$, or $0001$, or $0010$.
So the sequence $R^{(2)}_{00}$ is the union of the three sequences $3A+2\Id-5, 3A+2\Id-4$, and $3A+2\Id-3$.

\subsection{Densities}\label{sec:dens}

For a natural number $m$ let $w$ be a word of $0$'s and $1$'s, containing no $11$. By Theorem \ref{th:Zeck-w} we know that the sequence $R_w$ of occurrences of numbers $N$ such that  the $m$ lowest digits of the Zeckendorf expansion of $N$ are equal to $w$, is a  generalized Beatty sequence  given by
$$ R_w=F_{m+w_{m-1}}A+F_{m-1+w_{m-1}}\Id+\gamma_w, $$
for some negative integer $\gamma_w$.

As a result of this,  any word $w$ of length $m$ has a density $f_w$ of occurrence in the set of natural numbers, for respectively $w_{m-1}=0$, and $w_{m-1}=1$:
\begin{equation}\label{eq:freq0}
f_w=\frac1{F_{m}\varphi+F_{m-1}}=\varphi^{-m}\quad \text{or } \quad f_w=\frac1{F_{m+1}\varphi+F_{m}}=\varphi^{-m-1}.
\end{equation}
Here the equality
\begin{equation}\label{eq:FibPhi}
\varphi^{m}=F_{m}\varphi+F_{m-1}
\end{equation}
is easily proved by induction.




\begin{proposition}\label{prop:dens-gen} For any natural number $m$, fix a word $w=w_{m-1}\dots w_0$ of $0$'s and $1$'s, containing no $11$.
Let $k\ge 0$ be an integer. Let $f^{(k)}_w$ be the density of $R^{(k)}_w$ in $\mathbb{N}$. Then
$$ f^{(k)}_w=F_{k+2-w_0}\varphi^{-k-m-w_{m-1}}.$$
\end{proposition}

\noindent {\it Proof:} For $k=0$ this is Equation (\ref{eq:freq0}), since $F_1=F_2=1$. For $k>0$, this follows directly from Theorem \ref{th:gen-pos}, with Equation (\ref{eq:FibPhi}). \hfill $\Box$

\medskip

\noindent {\bf Remark:} There are $F_m$ words $w=0\dots0$, $F_{m-1}$ words $w=0\dots1$ or $w=1\dots0$, and $F_{m-2}$ words $w=1\dots1$.
So the total density of all words of length $m$ is equal to
$$F_mF_{k+2}\varphi^{-k-m}+F_{m-1}F_{k+1}\varphi^{-k-m}+F_{m-1}F_{k+2}\varphi^{-k-m-1}+F_{m-2}F_{k+1}\varphi^{-k-m-1}.$$
That this is equal to 1 follows from three instances of Equation (\ref{eq:FibPhi}), and from the well-known Fibonacci number relation\;
$F_mF_n+F_{m+1}F_{n+1}=F_{m+n+1}.$

\noindent AMS Classification Numbers: 11D85, 68R10, 11A63, 11B39

\end{document}